\documentclass[12pt]{amsart}

\usepackage{graphicx}
\usepackage{amsfonts}
\usepackage{hyperref}
\usepackage{textcomp}
\usepackage{amssymb}
\usepackage{color}
\usepackage{graphicx}

\newtheorem{theorem}{Theorem}

\newtheorem{proposition}[theorem]{Proposition}

\date{}
\title{Permutations containing many patterns}
\author[Albert]{M.~H.~Albert}
\address{Department of Computer Science \\%
University of Otago} 
\email{malbert@cs.otago.ac.nz}
\author[Coleman]{Micah Coleman}
\address{Department of Mathematics \\%
University of Florida} 
\author[Flynn]{Ryan Flynn}
\address{Department of Mathematics \\%
University of Florida}  
\author[Leader]{Imre Leader}
\address{Department of Pure Mathematics and Mathematical Statistics \\%
University of Cambridge}
\begin{document}
\begin{abstract}
It is shown that the maximum number of patterns that can occur in a permutation of length $n$ is asymptotically $2^n$.
This significantly improves a previous result of Coleman.
\end{abstract}
\maketitle

\section{Introduction}

Given a sequence $\mathbf{t} = t_1, t_2, \ldots, t_k$ of distinct elements from some totally ordered set, there is a unique permutation $\tau$ of $[k] = \{1, 2, \ldots , k \}$ with the property that for all $1 \leq i, j \leq k$, $t_i < t_j$ if and only if $\tau(i) < \tau(j)$. We call $\tau$ the \emph{pattern} of $\mathbf{t}$. For example, the pattern of $5, \, 10, \, 2$ written in one line notation is $231$. In other words, the sequence representing $\tau$ is obtained from $\mathbf{t}$ simply by replacing each element of $\mathbf{t}$ by its rank in $\mathbf{t}$.

Let $\sigma$ be a permutation of length $n$, written in one-line notation as $\sigma(1) \sigma(2) \cdots \sigma(n)$, and thought of as a sequence of length $n$. For each non-empty subset $X$ of $[n]$ define $\sigma_X$ to be the pattern of that subsequence of $\sigma$ whose indices belong to $X$. Define:
\[
P (\sigma) =   \{ \sigma_X \, : \, \emptyset \neq X \subseteq [n] \} .
\]
That is, $P(\sigma)$ is the set of patterns that occur in $\sigma$. Also define $h(n)$ to be the maximum value of $|P (\sigma)|$ taken over all permutations $\sigma$ of length $n$. 

Trivially, $h(n) \leq 2^n -1$. Slightly more precisely, for any permutation $\sigma$ of length $n$:
\begin{equation} \label{mainInequality}
\left| P (\sigma) \right| \leq \sum_{k = 1}^n \min \left(k!, \binom{n}{k} \right)
\end{equation}
since not more than $k!$ patterns of length $k$ can occur. However, the expression on the right hand side of this inequality is easily seen to be asymptotically $2^n$. At the 2003 conference on Permutation Patterns, Herb Wilf raised the issue of determining the (asymptotic) behaviour of $h(n)$, and exhibited a sequence of permutations which established that $h(n)$ exceeded the $n^{\rm{th}}$ Fibonacci number. Micah Coleman then demonstrated in \cite{Coleman:04} a sequence of permutations $\pi_n$, for $n$ a perfect square,\footnote{We have adjusted the notation slightly from that of \cite{Coleman:04} --- what was there called $\pi_k$ we are calling $\pi_{k^2}$ so that the subscript is equal to the length of the permutation.} for which:
\[
| P (\pi_n) | > 2^{n - 2 \sqrt{n} + 1}.
\]
Of course this establishes that $h(n)^{1/n} \rightarrow 2$ (for all $n$, not just perfect squares, using the fact that $h(n)$ is non decreasing). However, this left open the question of whether or not $h(n)/2^n$ tends to $1$ as $n$ tends to infinity.

In this paper, we refine the counting arguments concerning the number of patterns in $\pi_n$, for $n$ an even perfect square, and then extend the construction to all other values of $n$, in order to show that $| P (\pi_n) |/2^n \rightarrow 1$.  Indeed, we will obtain:
\[
h(n) > 2^n \left( 1 - 6 \sqrt{n} \, 2^{-\sqrt{n}/2}  \right)
\]
for all positive integers $n$.

\section{The main construction}

\begin{figure}[ht]
\begin{center}
\begin{tabular}{cc}
\includegraphics[width=5cm]{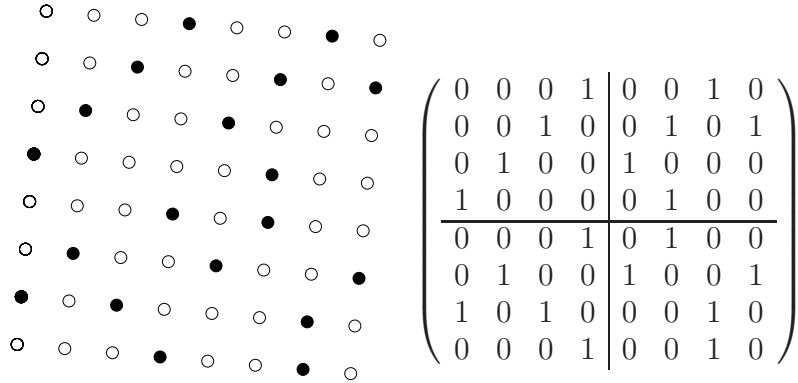} &
\raisebox{60pt}{\begin{minipage}[c]{5cm}
\[
\left(
\begin{array}{cccc|cccc}
0 & 0 & 0 & 1 & 0 & 0 & 1 & 0 \\
0 & 0 & 1 & 0 & 0 & 1 & 0 & 1 \\
0 & 1 & 0 & 0 & 1 & 0 & 0 & 0 \\
1 & 0 & 0 & 0 & 0 & 1 & 0 & 0 \\ \hline
0 & 0 & 0 & 1 & 0 & 1 & 0 & 0 \\
0 & 1 & 0 & 0 & 1 & 0 & 0 & 1 \\
1 & 0 & 1 & 0 & 0 & 0 & 1 & 0 \\
0 & 0 & 0 & 1 & 0 & 0 & 1 & 0 
\end{array}
\right)
\]
\end{minipage}}
\end{tabular}
\end{center}
\caption{The graph of the permutation $\pi_{64}$, an ample subset of its elements indicated by filled circles, together with the corresponding matrix divided into its corner submatrices.}
\label{fig01}
\end{figure}

Let $k$ be a positive integer and let $n = 4k^2$. Let $s$ be the sequence:
\[
s = (2k) \: (4k) \: (6k) \: \cdots \: (4k^2) 
\]
and consider the permutation $\pi_n$ which in one line notation is defined by:
\[
\pi_n = s \: (s-1) \: (s-2) \: \cdots \: (s - 2k + 1).
\]
Here $s-i$ indicates the sequence obtained by subtracting $i$ from each element of $s$. Generally, we will suppress the subscript on $\pi_n$ when there is no risk of confusion. Informally, the graph of $\pi$ is obtained by taking a standard orthogonal $2k \times 2k$ grid and rotating it slightly in the clockwise direction around its lower left hand corner. We associate to each subset $X$ of (the indices of) $\pi$ a $2k \times 2k$ 0-1 matrix, $M_X$, whose $1$ entries correspond to the elements of the subset. We also view $M_X$ as being partitioned into four $k \times k$ submatrices (called the \emph{corner submatrices}) in the usual way, that is, so that they form a $2 \times 2$ block decomposition of $M_X$. We say that $X$ (or $M_X$) is \emph{ample} if each $k \times k$ corner submatrix of $M_X$ has no zero rows or zero columns. An example is shown in Figure~\ref{fig01}.

\begin{proposition}
The number of ample matrices is greater than 
\[
2^n \left(1 - \frac{4\sqrt{n}}{2^{\sqrt{n}/2}} \right).
\]
\end{proposition}

\begin{proof}
Recall that $n = 4k^2$. Suppose that we sample an $n \times n$ 0-1 matrix uniformly at random from among all $n \times n$ 0-1 matrices. The probability that any particular row or column sum of one of the corner submatrices is 0 is $1/2^k$. There are $8k$ such sums which must all be non zero in order for the matrix to be ample. However, the probability that any of them are 0 is less than $8k/2^k$. So, the probability that all are non zero is greater than
\[
1 - \frac{8k}{2^k} = 1 - \frac{4\sqrt{n}}{2^{\sqrt{n}/2}},
\]
which is equivalent to the stated result.
\end{proof}

\begin{proposition}
Let $X$ and $Y$ be ample sets. Then $\pi_X = \pi_Y$ implies $X = Y$.
\end{proposition}

\begin{proof}
We must show that, if $X$ is ample, then it can be reconstructed from just the permutation $\pi_X$. Since $X$ is ample, the column sum of both the top half and bottom half of each column of $M_X$ is non zero. Therefore, there are $2k-1$ descents in $\pi_X$, corresponding to the transitions between columns of $M_X$. Thus, we can associate the elements of $\pi_X$ with their correct columns. However, this argument applies equally well to the rows of $M_X$ --- as is most easily seen by considering $\pi^{-1}$. Determining the row and column that represents each element of $\pi_X$ is exactly the same as reconstructing $X$.
\end{proof}

Combining these two results we have:

\begin{theorem}
\label{mainTheorem}
If $n$ is an even perfect square, then
\[
h(n) > 2^n \left( 1 - \frac{4\sqrt{n}}{2^{\sqrt{n}/2}} \right).
\]
\end{theorem}

We will refer to the second term inside the parentheses above as the \emph{correction term} for this estimate.

\section{Refinements}

It is easy to extend the above arguments to give lower bounds on
$h(n)$ that are valid for \emph{all} values of $n$. We can do this by using the basic construction of the previous section, and adding some extra elements in appropriate places to construct permutations $\pi_n$ of length $n$ that contain many patterns. 

\begin{figure}[ht]
\begin{center}
\includegraphics[width=5cm]{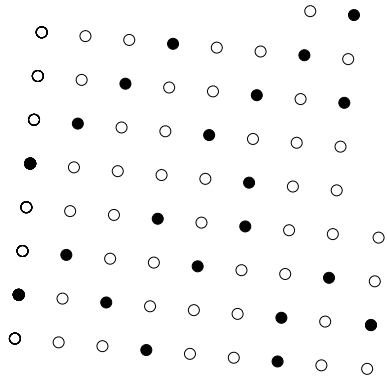}
\end{center}
\caption{The graph of the permutation $\pi_{70}$, together with the matrix associated with a particular ample subset of its elements indicated by filled circles.}
\label{fig02}
\end{figure}

First suppose that $n = 4k^2 + l$ where $0 < l < 2k$. Take the grid associated to the permutation $\pi_{4k^2}$ and add a (partial) column on the right hand side at the bottom containing not more than $k$ elements, and, if necessary, a partial row on top at the right hand side, also not containing more than $k$ elements, so that the total number of elements added is $l$. As before, rotate this grid slightly, and view the result as the graph of a permutation, $\pi_n$. An example is shown in Figure~\ref{fig02}. Call the elements of this permutation arising from the original grid defining $\pi_{4k^2}$ the \emph{main} elements, and the remaining elements the \emph{extra} elements. Define a subset of the indices of $\pi_n$ to be \emph{ample} if its intersection with the main elements would be ample for $\pi_{4k^2}$.

\begin{proposition}
Let $X$ and $Y$ be ample sets. Then $\pi_X = \pi_Y$ implies $X = Y$.
\end{proposition}

\begin{proof}
As before, we must describe how to reconstruct $X$ from $\pi_X$. However, we can identify the extra elements (and hence the main elements) in $\pi_X$. If there are any belonging to the new partial column, then they are exactly the elements following the $(2k)^{\rm{th}}$ descent, while those belonging to the new partial row, if such exist, are exactly those lying above the maximum element of the first $k$ columns. Since the main elements form an ample subset of $\pi_{4k^2}$ we can use the preceding result to identify their values. Once the values of the main elements are known, so are the values of the extra elements.
\end{proof}

Therefore, for such $n$,
\[
h(n) \geq |P (\pi_n) | > 2^{4k^2} \left( 1 - \frac{8k}{2^k} \right) 2^l .
\]
Certainly $k \leq \sqrt{n}/2$, but also $(2k+1/2)^2 > n$ so $ k > (\sqrt{n} - 1/2)/2$. Applying these estimates we obtain:
\[
h(n) > 2^n \left( 1 - \frac{2^{9/4} \sqrt{n}}{2^{\sqrt{n}/2}} \right).
\]
This differs from our previous estimate by a factor of $2^{1/4}$ in the correction term.

For $n = 4k^2 + 2k$, we switch to a grid consisting of $2k+1$ columns of size $2k$ and define $\pi_n$ appropriately. As in the previous section, we define the four corner submatrices, except now those on the right hand side of the matrix are $k \times (k+1)$ instead of $k \times k$. The probability of a subset of the matrix not being ample is not as much as:
\[
\frac{2(2k+1)}{2^k} + \frac{2k}{2^k} + \frac{2k}{2^{k+1}} = \frac{7k+2}{2^k}.
\]
Using the same bounds as before (which still apply) plus trivial estimates for $k \leq 2$ it is easy to check that the bound
\[
h(n) > 2^n \left( 1 - \frac{2^{9/4} \sqrt{n}}{2^{\sqrt{n}/2}} \right)
\]
still applies in this case. We can proceed from this point with the half-row/half-column construction again (possibly at a penalty of another factor of $2^{1/4}$ in the correction term) as far as $n = (2k+1)^2$. At this point we pause for a detailed re-evaluation. In a $(2k+1) \times (2k+1)$ grid, divided into corner submatrices of sizes $k \times k$, $k \times (k+1)$, $(k+1) \times k$ and $(k+1) \times (k+1)$, the probability that a subset is not ample is less than:
\[
\frac{2k}{2^k} + 2 \left( \frac{k}{2^{k+1}} + \frac{k+1}{2^k} \right) + \frac{2(k+1)}{2^{k+1}} = \frac{6k+3}{2^k}.
\]
Since $k = (\sqrt{n}-1)/2$, this equals
\[
\frac{(3 \sqrt{2}) \sqrt{n}}{2^{\sqrt{n}/2}}.
\]
We can pursue these constructions through to the next even perfect square, and, allowing  for a further penalty of $\sqrt{2}$ in the correction term (which we leave to the reader to verify is generous), obtain:
\begin{theorem}
For all positive integers $n$,
\[
h(n) > 2^n \left( 1 - \frac{6 \sqrt{n}}{2^{\sqrt{n}/2}} \right).
\]
\end{theorem}

\section{Conclusions}

It would be interesting to know just how close to $2^n$ the value of
$h(n)$ actually is. A more careful analysis of the various steps in moving from one square grid to the next might well provide a small improvement in the constant factor of the correction term of our estimate. Similarly, an analysis of conditions weaker than ample which none the less would allow for a reconstruction result might actually improve the asymptotic form of the correction term. However, the simplicity of the main construction (for $n = 4k^2$) and of the proof that ample subsets can be reconstructed from their patterns, together with the lack of any great \emph{need} for more precise estimates of $h(n)$ somewhat dampens our enthusiasm for further investigations in that direction. Of perhaps greater interest would be to investigate the distribution of the statistic $| P (\pi) |$ as $\pi$ ranges over permutations of length $n$. 

We would like to thank Herb Wilf for having posed such an interesting problem!

\bibliographystyle{plain}
\bibliography{maxPatterns}

\begin{thebibliography}{1}

\bibitem{Coleman:04}
Micah Coleman.
\newblock An answer to a question by {W}ilf on packing distinct patterns in a
  permutation.
\newblock {\em Electron. J. Combin.}, 11(1):Note 8, 4 pp. (electronic), 2004.

\end{thebibliography}

\end{document}